\documentclass[12pt]{amsart}

\newtheorem{theorem}{Theorem}[section]

\theoremstyle{definition}
\newtheorem{definition}{Definition}[section]

\theoremstyle{remark}

\numberwithin{equation}{section}

\begin{document}
\title[Ricci Solitons on submanifolds of $(LCS)_n$-Manifolds]{Ricci Solitons on submanifolds of $(LCS)_n$-Manifolds}
\author[S. K. Hui, R. Prasad and T. Pal]{SHYAMAL KUMAR HUI, RAJENDRA PRASAD and TANUMOY PAL}
.
\begin{abstract}
The present paper deals with the study of Ricci solitons on invariant and anti-invariant
submanifolds of $(LCS)_n$-manifolds with respect to Riemannian connection as well as quarter symmetric metric connection.
\end{abstract}
\maketitle
\section{Introduction}
\indent In 1982, Hamilton \cite{HAM1} introduced the notion of Ricci flow to find a canonical metric on a smooth manifold. Then Ricci flow
has become a powerful tool for the study of Riemannian manifolds, especially for manifolds with positive curvature.
Perelman \cite{PER1} used Ricci flow and its surgery to prove Poincare conjecture.
The Ricci flow is an evolution equation for metrics on a Riemannian manifold defined as follows:
\begin{equation*}
\frac{\partial}{\partial t}g_{ij}(t) = - 2 R_{ij}.
\end{equation*}
\indent A Ricci soliton emerges as the limit of the solutions of the Ricci flow. A solution to the Ricci flow is called
Ricci soliton if it moves only by a one parameter group of diffeomorphism and scaling. A Ricci soliton $(g,V,\lambda)$
on a Riemannian manifold $(M,g)$ is a generalization of an Einstein metric such that \cite{HAM2}
\begin{equation}
\label{eqn1.3}
\pounds_{V}g + 2S + 2\lambda g = 0,
\end{equation}
where $S$ is the Ricci tensor, $\pounds_{V}$ is the Lie derivative operator along the vector field $V$
on $M$ and $\lambda$ is a real number. The Ricci soliton is said to be shrinking, steady and
expanding according as $\lambda$ is negative, zero and positive respectively.\\
\indent During the last two decades, the geometry of Ricci solitons has been
the focus of attention of many mathematicians. In particular, it has
become more important after Perelman applied Ricci solitons to solve
the long standing Poincare conjecture posed in 1904. In \cite{SHAR} Sharma studied the Ricci solitons
in contact geometry. Thereafter, Ricci solitons in
contact metric manifolds have been studied by various authors such as Bejan and Crasmareanu \cite{BEJA},
Hui et. al (\cite{CHS},\cite{SKH2}-\cite{SKH6}, \cite{SKH9}), Chen and Deshmukh \cite{CD}, Deshmukh et. al \cite{Detal},
He and Zhu \cite{HE}, Tripathi \cite{TRIP} and many others.\\
\indent \indent In 2003, Shaikh \cite{SHAIKH2} introduced the notion of Lorentzian
concircular structure manifolds (briefly, $(LCS)_n$-manifolds),
with an example, which generalizes the notion of LP-Sasakian
manifolds introduced by Matsumoto \cite{8} and also by Mihai and
Rosca \cite{9}. Then Shaikh and Baishya (\cite{SHAIKH3}, \cite{SHAIKH4})
investigated the applications of $(LCS)_n$-manifolds to the general
theory of relativity and cosmology. The $(LCS)_n$-manifolds is also
studied by Hui \cite{SKH}, Hui and Atceken \cite{HUIAT}, Shaikh and his
co-authors (\cite{SHAIKH1}-\cite{SHAIKH9}) and many others.\\
\indent In modern analysis, the geometry of submanifolds has become a subject of growing interest
for its significant applications in applied mathematics and theoretical physics.
The present paper deals with the study of Ricci solitons on submanifolds of $(LCS)_n$-manifolds.
The paper is organised as follows. Section 2 is concerned with some preliminaries. Section 3 is devoted
to the study of Ricci solitons on invariant and anti-invariant submanifolds of $(LCS)_n$-manifolds.\\
\indent In 1924, Friedman and Schouten \cite{FRSC} introduced the notion of semi-symmetric linear
connection on a differentiable manifold. In 1932, Hayden \cite{HAYDEN} introduced the idea of metric
connection with torsion on a Riemannian manifold. In 1970, Yano \cite{YANO} studied some curvature
tensors and conditions for semi-symmetric connections in Riemannian manifolds. In 1975, Golab \cite{GOLAB}
defined and studied quarter symmetric linear connection on a differentiable manifold. A linear
connection $\overline{\nabla}$ in an n-dimensional Riemannian manifold is said to be a quarter
symmetric connection \cite{GOLAB} if torsion tensor $T$ is of the form
\begin{equation}\label{1.1}
T(X,Y)= \overline{\nabla}_XY - \overline{\nabla}_YX - [X,Y] = A(Y) K(X)-A(X)K(Y)
\end{equation}
where $A$ is a 1-form and $K$ is a tensor of type (1,1). If a quarter symmetric linear connection $\overline{\nabla}$ satisfies the condition
\begin{equation*}
(\overline{\nabla}_Xg)(Y,Z)=0
\end{equation*}
for all $X$, $Y$, $Z\in \chi(M)$, where $\chi(M)$ is a Lie algebra of vector fields on the
manifold $M$, then $\overline{\nabla}$ is said to be a quarter symmetric metric connection.
For a contact metric manifold admitting quarter symmetric connection, we can
take $A=\eta$ and $K=\phi$ and hence  (\ref{1.1}) takes in the form:
\begin{equation}\label{1.2}
T(X,Y)=\eta(Y)\phi X-\eta(X)\phi Y.
\end{equation}
The relation between Levi-Civita connection $\nabla$ and quarter symmetric metric connection $\overline{\nabla}$
of a contact metric manifold is given by
\begin{equation}
\label{1.3}
\overline{\nabla}_XY=\nabla_XY-\eta(X)\phi Y.
\end{equation}
\indent Recently Hui, Piscoran and Pal \cite{HPP} studied invariant submanifolds of $(LCS)_n$-manifolds with
respect to quarter symmetric metric connection. Ricci solitons on invariant and anti-invariant
submanifolds of $(LCS)_n$-manifolds with respect to quarter symmetric metric connections studied in section 4 of the paper.
\section{preliminaries}
\par An $n$-dimensional Lorentzian manifold $\widetilde{M}$ is a smooth connected
paracompact Hausdorff manifold with a Lorentzian metric $g$, that
is, $\widetilde{M}$ admits a  smooth symmetric tensor field $g$ of type (0,2)
such that for each point $p\in \widetilde{M}$, the tensor $g_{p}:T_{p}\widetilde{M}\times
T_{p}\widetilde{M}$ $\rightarrow\mathbb{R}$ is a non-degenerate inner product of
signature $(-,+,\cdots,+)$, where $T_{p}\widetilde{M}$ denotes the tangent
vector space of $\widetilde{M}$ at $p$ and $\mathbb{R}$ is the real number
space. A non-zero vector $v$ $\in T_{p}\widetilde{M}$ is said to be timelike
(resp., non-spacelike, null, spacelike) if it satisfies $g_{p}(v,v)
< 0$ (resp, $\leq $ 0, = 0, $> 0$) \cite{NIL}.
\begin{definition}
In a Lorentzian manifold $(\widetilde{M},g)$ a vector field $P$ defined by
\begin{equation*}
g(X,P)=A(X)
\end{equation*}
for any $X\in\Gamma(T\widetilde{M})$, is said to be a concircular vector field \cite{15} if
\begin{equation*}
(\widetilde{\nabla}_{X}A)(Y)=\alpha \{g(X,Y)+\omega(X)A(Y)\},
\end{equation*}
where $\alpha$ is a non-zero scalar and $\omega$ is a closed 1-form
and $\widetilde{\nabla}$ denotes the operator of covariant
differentiation with respect to the Lorentzian metric $g$.
\end{definition}
Let $\widetilde{M}$ be an $n$-dimensional Lorentzian manifold admitting a unit
timelike concircular vector field $\xi$, called the characteristic
vector field of the manifold. Then we have
\begin{equation}
\label{2.1}
g(\xi, \xi)=-1.
\end{equation}
Since $\xi$ is a unit concircular vector field, it follows that
there exists a non-zero 1-form $\eta$ such that for
\begin{equation}
\label{2.2}
g(X,\xi)=\eta(X),
\end{equation}
the equation of the following form holds
\begin{equation}
\label{2.3}
(\widetilde\nabla _{X}\eta)(Y)=\alpha \{g(X,Y)+\eta(X)\eta(Y)\},
\ \ \ (\alpha\neq 0)
\end{equation}
\begin{equation}
\label{2.4}
\widetilde\nabla _{X}\xi = \alpha \{X +\eta(X)\xi\}, \ \ \ \alpha\neq 0,
\end{equation}
for all vector fields $X$, $Y$, where $\widetilde{\nabla}$ denotes the
operator of covariant differentiation with respect to the Lorentzian
metric $g$ and $\alpha$ is a non-zero scalar function satisfies
\begin{equation}
\label{2.5}
{\widetilde\nabla}_{X}\alpha = (X\alpha) = d\alpha(X) = \rho\eta(X),
\end{equation}
$\rho$ being a certain scalar function given by $\rho=-(\xi\alpha)$.
Let us take
\begin{equation}
\label{2.6}
\phi X=\frac{1}{\alpha}\widetilde\nabla_{X}\xi,
\end{equation}
then from (\ref{2.4}) and (\ref{2.6}) we have
\begin{equation}
\label{2.7} \phi X = X+\eta(X)\xi,
\end{equation}
\begin{equation}
\label{2.8}
g(\phi X,Y) = g(X,\phi Y),
\end{equation}
from which it follows that $\phi$ is a symmetric (1,1) tensor and
called the structure tensor of the manifold. Thus the Lorentzian
manifold $\widetilde{M}$ together with the unit timelike concircular vector
field $\xi$, its associated 1-form $\eta$ and an (1,1) tensor field
$\phi$ is said to be a Lorentzian concircular structure manifold
(briefly, $(LCS)_{n}$-manifold), \cite{SHAIKH2}. Especially, if we take
$\alpha=1$, then we can obtain the LP-Sasakian structure of
Matsumoto \cite{8}. In a $(LCS)_{n}$-manifold $(n>2)$, the following
relations hold \cite{SHAIKH2}:
\begin{equation}
\label{2.9}
\eta(\xi)=-1,\ \ \phi \xi=0,\ \ \ \eta(\phi X)=0,\ \ \
g(\phi X, \phi Y)= g(X,Y)+\eta(X)\eta(Y),
\end{equation}
\begin{equation}
\label{2.10}
\phi^2 X= X+\eta(X)\xi,
\end{equation}
\begin{equation}
\label{2.11}
\widetilde{S}(X,\xi)=(n-1)(\alpha^{2}-\rho)\eta(X),
\end{equation}
\begin{equation}
\label{2.12}
\widetilde{R}(X,Y)\xi=(\alpha^{2}-\rho)[\eta(Y)X-\eta(X)Y],
\end{equation}
\begin{equation}
\label{2.13}
\widetilde{R}(\xi,Y)Z=(\alpha^{2}-\rho)[g(Y,Z)\xi-\eta(Z)Y],
\end{equation}
\begin{equation}
\label{2.14}
(\widetilde{\nabla}_{X}\phi)Y=\alpha\{g(X,Y)\xi+2\eta(X)\eta(Y)\xi+\eta(Y)X\},
\end{equation}
\begin{equation}
\label{2.15}
(X\rho)=d\rho(X)=\beta\eta(X),
\end{equation}
\begin{equation}
\label{2.16}
\widetilde{R}(X,Y)Z =\phi \widetilde{R}(X,Y)Z +(\alpha^{2}-\rho)\{g(Y,Z)\eta(X)-g(X,Z)\eta(Y)\}\xi
\end{equation}
for all $X,\ Y,\ Z\in\Gamma(T\widetilde{M})$ and $\beta = -(\xi\rho)$ is a scalar function,
where $\widetilde{R}$ is the curvature tensor and $\widetilde{S}$ is the Ricci tensor of the manifold.\\
\indent Let $M$ be a submanifold of dimension $m$ of a $(LCS)_n$-manifold $\widetilde{M}$ $(m<n)$ with induced
metric $g$. Also let $\nabla$ and $\nabla^{\perp}$ be the induced
connection on the tangent bundle $TM$ and the normal bundle
$T^{\perp}M$ of $M$ respectively. Then the Gauss and Weingarten
formulae are given by
\begin{equation}\label{2.17}
\widetilde{\nabla}_{X}Y = \nabla_{X}Y + h(X,Y)
\end{equation}
and
\begin{equation}\label{2.18}
\widetilde{\nabla}_{X}V = -A_{V}X + \nabla^{\perp}_{X}V
\end{equation}
for all $X,\ Y \in\Gamma(TM)$ and $V\in\Gamma(T^{\perp}M)$, where $h$
and $A_V$ are second fundamental form and the shape operator
(corresponding to the normal vector field $V$) respectively for the
immersion of $M$ into $\widetilde{M}$. The second fundamental form $h$ and the
shape operator $A_V$ are related by \cite{YANO3}
\begin{equation}\label{2.19}
g(h(X,Y),V) = g(A_{V}X,Y),
\end{equation}
for any $X,\ Y \in\Gamma(TM)$ and $V\in\Gamma(T^{\perp}M)$. We note that $h(X,Y)$ is bilinear and since
$\nabla_{fX}Y = f\nabla_{X}Y$ for any smooth function $f$ on a manifold, we have
\begin{equation}\label{2.20}
h(fX,Y) = f h(X,Y).
\end{equation}
\noindent The mean curvature vector $H$ on $M$ is given by $H=\frac{1}{m}\displaystyle\sum_{i=1}^{m}g(e_i,e_i)$,
where $\{e_1,e_2,\cdots,e_m\}$ is a local orthonormal frame of vector fields on $M$.\\
\noindent A submanifold $M$ of a $(LCS)_n$-manifold $\widetilde{M}$ is said to be totally umbilical if
\begin{equation}\label{2.21}
h(X,Y)=g(X,Y)H
\end{equation}
for any vector fields $X$, $Y$ $\in$ $TM$. Moreover if $h(X,Y)=0$ for all $X,\
Y\ \in TM$, then $M$ is said to be totally geodesic and if $H=0$ then $M$ is minimal in $\widetilde{M}$.

\indent Analogous to almost Hermitian manifolds, the invariant and anti-invariant submanifols are depend on
the behaviour of almost contact metric structure $\phi$.\\
\indent A submanifold $M$ of an almost contact metric manifold $\widetilde{M}$ is said to be invariant if
the structure vector field $\xi$ is tangent to $M$ at every point of $M$ and $\phi X$ is tangent to $M$ for
every vector field $X$ tangent to $M$ at evey point of $M$. i.e.
$ \phi(TM)\subset TM$ at evey point of $M$.\\
\indent On the other hand, $M$ is said to be anti-invariant if for any $X$ tangent to $M$,
$\phi X$ is normal to $M$, i.e., $\phi TM\subset T^\bot M$ at every point of $M$,
where $T^\bot M$ is the normal bundle of $M$.\\
\indent Let $\overline{\widetilde{\nabla}}$ be a linear connection and $\widetilde{\nabla}$ be the Levi-Civita connection of a $(LCS)_n$-manifold
$\widetilde{M}$ such that
\begin{equation}\label{2.22}
\overline{\widetilde{\nabla}}_XY=\widetilde{\nabla}_XY+U(X,Y),
\end{equation}
where $U$ is a (1,1) type tensor and $X,\ Y\in \Gamma(T\widetilde{M})$.
For $\overline{\widetilde{\nabla}}$ to be a quarter symmetric metric connection on $\widetilde{M}$, we have
\begin{equation}\label{2.23}
U(X,Y)=\frac{1}{2}[T(X,Y)+T^\prime(X,Y)+ T^\prime(Y,X)],
\end{equation}
where
\begin{equation}\label{2.24}
g(T^\prime(X,Y),Z)=g(T(Z,X),Y).
\end{equation}
From (\ref{1.2}) and (\ref{2.24}) we get
\begin{equation}\label{2.26}
T^\prime(X,Y)=\eta(X)\phi Y-g(Y,\phi X)\xi.
\end{equation}
So,
\begin{equation}\label{2.27}
U(X,Y)=\eta(Y)\phi X-g(Y,\phi X)\xi.
\end{equation}
Therefore a quarter symmetric metric connection $\overline{\widetilde{\nabla}}$ in a $(LCS)_n$-manifold $\widetilde{M}$ is given by
\begin{equation}\label{2.28}
\overline{\widetilde{\nabla}}_XY=\widetilde{\nabla}_XY+ \eta(Y)\phi X-g(\phi X,Y)\xi.
\end{equation}
\par Let $\overline{\widetilde{R}}$ and $\widetilde{R}$ be the curvature tensors of a $(LCS)_n$-manifold $\widetilde{M}$
with respect to the quarter symmetric metric connection $\overline{\widetilde{\nabla}}$ and
the Levi-Civita connection $\widetilde{\nabla}$ respectively. Then we have
\begin{eqnarray}\label{2.29}
  \overline{\widetilde{R}}(X,Y)Z&=& \widetilde{R}(X,Y)Z+(2\alpha-1)\left[g(\phi X,Z)\phi Y-g(\phi Y,Z)\phi X \right]  \\
 \nonumber &&+\alpha\left[\eta(Y)X-\eta(X)Y\right]\eta(Z) \\
\nonumber &&  +\alpha \left[ g(Y,Z)\eta(X)-g(X,Z)\eta(Y)\right]\xi,
\end{eqnarray}
where $\overline{\widetilde{R}}(X,Y)Z=\overline{\widetilde{\nabla}}_X \overline{\widetilde{\nabla}}_Y Z-\overline{\widetilde{\nabla}}_Y \overline{\widetilde{\nabla}}_XZ-\overline{\widetilde{\nabla}}_{[X,Y]}Z$
and $X,\ Y,\ Z\in \chi(\widetilde{M})$.
\par By suitable contraction we have from (\ref{2.29}) that
\begin{eqnarray}\label{2.30}
  \overline{\widetilde{S}}(Y,Z) &=& \widetilde{S}(Y,Z)+(\alpha-1)g(Y,Z)+(n\alpha-1)\eta(Y)\eta(Z) \\
  \nonumber&&-(2\alpha-1)a g(\phi Y,Z),
\end{eqnarray}
where $\overline{\widetilde{S}}$ and $\widetilde{S}$ are the Ricci tensors of $\widetilde{M}$ with respect to
 $\overline{\widetilde{\nabla}}$ and $\widetilde{\nabla}$ respectively and $a=trace\phi$.
\section{Ricci solitons on submanifolds of $(LCS)_n$-Manifolds}
Let us take $(g,\xi,\lambda)$ be a Ricci soliton on a submanifold $M$ of a $(LCS)_n$-manifold $\widetilde{M}$.\\
 Then we have
 \begin{equation}\label{3.1}
 (\pounds_\xi g)(Y,Z)+2S(Y,Z)+2\lambda g(Y,Z)=0.
 \end{equation}
 From (\ref{2.6}) and (\ref{2.17}) we get
 \begin{equation}\label{3.2}
 \alpha \phi X= \widetilde{\nabla}_X\xi=\nabla_X\xi+h(X,\xi).
 \end{equation}
 If $M$ is invariant in $\widetilde{M}$, then $\phi X, \xi\in TM$ and therefore equating tangential and normal components of $(\ref{3.2})$ we get
 \begin{equation}\label{3.3}
 \nabla_X\xi=\alpha \phi X\ \  \text{and}\ \  h(X,\xi)=0.
 \end{equation}
 From (\ref{2.1}), (\ref{2.2}), (\ref{2.7}) and (\ref{3.3}) we get
 \begin{eqnarray}
 \label{3.4}(\pounds_\xi g)(Y,Z)&=&g(\nabla_Y\xi,Z)+g(Y,\nabla_Z\xi)\\
 \nonumber&&=2\alpha [g(Y,Z)+\eta(Y)\eta(Z)].
 \end{eqnarray}
 In view of (\ref{3.4}), (\ref{3.1}) yields
 \begin{equation}\label{3.5}
 S(Y,Z)=-(\alpha+\lambda)g(Y,Z)-\alpha \eta (Y)\eta(Z),
 \end{equation}
 which implies that $M$ is $\eta$-Einstein. Also from (\ref{2.20}) and (\ref{3.3}) we get $\eta(X) H=0$, i.e., $H=0$, since $\eta(X)\neq 0$.\\
 Consequently $M$ is minimal in $\widetilde{M}$. Thus we can state the following:
\begin{theorem}
If $(g,\xi,\lambda)$ is a Ricci soliton on an invariant submanifold $M$ of a $(LCS)_n$-manifold $\widetilde{M}$,
then $M$ is $\eta$-Einstein and also $M$ is minimal in $\widetilde{M}$.
\end{theorem}
From (\ref{3.3}) and using the formula
\begin{equation*}
R(X,Y)\xi=\nabla_X\nabla_Y\xi-\nabla_Y\nabla_X\xi-\nabla_{[X,Y]}\xi
\end{equation*}
we get
\begin{equation*}
R(X,Y)\xi=(\alpha^2-\rho)[\eta(Y)X-\eta(X)Y]
\end{equation*}
from which it follows that
\begin{equation}\label{3.6}
S(X,\xi)=(m-1)(\alpha^2-\rho)\eta(X) \ \ \ \ \text{for all X}.
\end{equation}
Putting $Z=\xi$ in (\ref{3.5}) and using (\ref{2.2}) and (\ref{3.6}) we get $\lambda=-(m-1)(\alpha^2-\rho)$. This leads to the following:
\begin{theorem}
A Ricci soliton $(g,\xi,\lambda)$ on an invariant submanifold of a $(LCS)_n$-manifold is shrinking, steady and expanding according as
$\alpha^2-\rho < 0$, $\alpha^2-\rho=0$ and $\alpha^2-\rho >0$ respectively.
\end{theorem}
\indent Again, if $M$ is anti-invariant in $\widetilde{M}$, then for any $X\in TM$, $\phi(X)\in T^\bot M$ and hence from (\ref{3.2}) we get
$\nabla_X\xi=0$ and $h(X,\xi)=\alpha \phi (X)$. Then $$(\pounds_\xi g)(Y,Z)=g(\nabla_Y\xi,Z)+(Y,\nabla_Z\xi)=0,$$ which means that $\xi$ is a Killing vector field and consequently (\ref{3.1}) yields $$S(Y,Z)=-\lambda g(Y,Z),$$ which implies that $M$ is Einstein. Thus we can state the following:
\begin{theorem}
If $(g,\xi,\lambda)$ is a Ricci soliton on an anti-invariant submanifold $M$ of a $(LCS)_n$-manifold $\widetilde{M}$,
then $M$ is Einstein and $\xi$ is Killing vector field.
\end{theorem}
Also, from $\nabla_X\xi=0$ we get $R(X,Y)\xi=0$ and hence $S(Y,\xi)=0$.
Again, we have $S(Y,\xi)=-\lambda\eta(Y)$. Therefore $\lambda=0$ and hence the Ricci soliton $(g,\xi,\lambda)$ is always steady.
This leads to the following:
\begin{theorem}
A Ricci soliton $(g,\xi,\lambda)$ on an anti-invariant submanifold $M$ of a $(LCS)_n$-manifold $\widetilde{M}$ is always steady.
\end{theorem}
\section{Ricci solitons on submanifolds of $(LCS)_n$-Manifolds with respect to quarter symmetric metric connection}
We now consider $(g,\xi, \lambda)$ is a Ricci soliton on a submanifold $M$ of a $(LCS)_n$-manifold $\widetilde{M}$
with respect to quarter symmetric metric connection, where $\overline{\nabla}$ is the induced connection on $M$ from the connection $\overline{\widetilde{\nabla}}$. Then we have
\begin{equation}\label{3.7}
(\overline{\pounds}_\xi g)(Y,Z)+2\overline{S}(Y,Z)+2\lambda g(Y,Z)=0.
\end{equation}
\indent Let $\overline{h}$ be the second fundamental form of $\overline{M}$ with respect to induced connection $\overline{\nabla}$. Then we have
\begin{equation}\label{3.8}
\overline{\widetilde{\nabla}}_XY=\overline{\nabla}_XY+\overline{h}(X,Y)
\end{equation}
and hence by virtue of (\ref{2.17}) and (\ref{2.26}) we get
\begin{equation}\label{3.9}
\overline{\nabla}_XY+\overline{h}(X,Y)=\nabla_XY+h(X,Y)+\eta(Y)\phi X-g(\phi X,Y)\xi.
\end{equation}
If $M$ is invariant submanifold of $\widetilde{M}$ then $\phi X,\ \xi\in TM$ for any $X\in TM$
and therefore equating tangential part from (\ref{3.7}) we get
\begin{equation}\label{3.10}
\overline{\nabla}_XY=\nabla_XY+\eta(Y)\phi X-g(\phi X,Y)\xi,
\end{equation}
which means $M$ admits quarter symmetric metric connection.\\
Also from (\ref{3.10}) we get $\overline{\nabla}_X\xi=(\alpha-1)\phi X$ and hence
\begin{eqnarray}\label{3.11}
(\overline{\pounds}_\xi g)(Y,Z)&=&g(\overline{\nabla}_Y\xi,Z)+g(Y,\overline{\nabla}_Z\xi)\\
\nonumber&&=2(\alpha-1)[g(Y,Z)+\eta(Y)\eta(Z)].
\end{eqnarray}
If $\overline{R}$ be the curvature tensor of submanifold $M$ with respect to induced connection $\overline{\nabla}$ of $(LCS)_n$-manifold
$\widetilde{M}$ with respect to quarter symmetric metric connection $\overline{\widetilde{\nabla}}$. Then we have,
\begin{eqnarray}\label{4.7}
  \overline{R}(X,Y)Z&=& R(X,Y)Z+(2\alpha-1)\left[g(\phi X,Z)\phi Y-g(\phi Y,Z)\phi X \right]  \\
 \nonumber &&+\alpha\left[\eta(Y)X-\eta(X)Y\right]\eta(Z) \\
\nonumber &&  +\alpha \left[ g(Y,Z)\eta(X)-g(X,Z)\eta(Y)\right]\xi,
\end{eqnarray}
where $\overline{R}(X,Y)Z=\overline{\nabla}_X \overline{\nabla}_Y Z-\overline{\nabla}_Y \overline{\nabla}_XZ-\overline{\nabla}_{[X,Y]}Z$.
Taking suitable contraction of above equation, we get
\begin{equation}\label{3.12}
\overline{S}(Y,Z)=S(Y,Z)+[\alpha(1-2a)+a]g(Y,Z)+[\alpha(m-2a)+a-1]\eta(Y)\eta(Z)
\end{equation}
Using (\ref{3.11}) and (\ref{3.12}) in (\ref{3.7}), we get
\begin{equation}\label{3.13}
S(Y,Z)=[2\alpha(a-1)+1-a-\lambda]g(Y,Z)+[\alpha(2a-m-1)+2-a]\eta(Y)\eta(Z),
\end{equation}
which implies that $M$ is $\eta$-Einstein.
\begin{theorem}
Let $(g,\xi,\lambda)$ be a Ricci soliton on an invariant submaniold $M$ of a $(LCS)_n$-manifold $\widetilde{M}$
with respect to quarter symmetric metric connection $\overline{\widetilde{\nabla}}$.
Let $\overline{\nabla}$ be the induced connectiuon on $M$ from the connection $\overline{\widetilde{\nabla}}$.
Then $M$ is $\eta$-Einstein with respect to Levi-Civita connection.
\end{theorem}
\indent Again, if $M$ is an anti-invariant submanifold of $\widetilde{M}$ with respect to quarter symmetric metric connection,
 then from (\ref{3.9}) we get $\overline{\nabla}_X\xi =0$. Consequently we get
 \begin{equation}\label{3.14}
 (\overline{\pounds}_\xi g)(Y,Z)=0.
 \end{equation}
 In view of (\ref{3.14}), (\ref{3.7}) yields
 \begin{equation*}
   \overline{S}(Y,Z)=-\lambda g(Y,Z),
 \end{equation*}
 which implies that $M$ is $\eta$-Einstein with respect to Riemannian connection by virtue of (\ref{3.12}). Thus we can state the following:
 \begin{theorem}
Let $(g,\xi,\lambda)$ be a Ricci soliton on an anti-invariant submaniold $\widetilde{M}$ of a $(LCS)_n$-manifold $\widetilde{M}$
with respect to quarter symmetric metric connection $\overline{\widetilde{\nabla}}$. Then $M$ is $\eta$-Einstein with respect
to induced Riemannian connection.
 \end{theorem}
\section{conclusion}
In this paper, we have studied invariant submanifolds of $(LCS)_n$-manifold $\widetilde{M}$ whose metric are Ricci solitons.
From Theorem 3.1, Theorem 3.3, Theorem 4.1 and Theorem 4.2, we can state the following:\\
\medskip
\noindent{\bf Theorem 5.1.} Let $(g,\xi,\lambda)$ be a Ricci soliton on a submanifold $M$ of a $(LCS)_n$-manifold $\widetilde{M}$.
Then the following holds:
\begin{center}
\begin{tabular}{|c|c|c|}
\hline nature of submanifold & connection of $\widetilde{M}$ & $M$ \\
\hline invariant & Riemannian & $\eta$-Einstein\\
\hline anti-invariant & Riemannian & Einstein\\
\hline invariant & quarter symmetric metric & $\eta$-Einstein\\
\hline anti-invariant & quarter symmetric metric & $\eta$-Einstein\\
\hline
\end{tabular}
\end{center}
\vspace{0.1in}
\noindent{\bf Acknowledgement:} The first author (S. K. Hui) gratefully acknowledges to
the SERB (Project No.: EMR/2015/002302), Govt. of India for financial assistance of the work.

\vspace{0.1in}
\noindent S. K. Hui\\
Department of Mathematics, The University of Burdwan, Golapbag, Burdwan -- 713104, West Bengal, India\\
E-mail: shyamal\_hui@yahoo.co.in, skhui@math.buruniv.ac.in\\

\vspace{0.1in}
\noindent R. Prasad\\
Department of Mathematics and Astronomy, University of Lucknow, Lucknow -- 226007, India.\\
E-mail: rp.manpur@rediffmail.com\\

\vspace{0.1in}
\noindent T. Pal\\
A. M. J. High School, Mankhamar, Bankura -- 722144, West Bengal, India\\
E-mail: tanumoypalmath@gmail.com
\end{document}